\documentclass[12pt,reqno]{amsart}

\usepackage{amsmath, amssymb, amscd, amsthm}
\usepackage[mathscr]{eucal}

\newdimen\unit\newdimen\psep\newcount\nd\newcount\ndx\newbox\dotb\newbox\ptbox
\newdimen\dx\newdimen\dy\newdimen\dxx\newdimen\dyy\newdimen\hgt
\newdimen\xoff\newdimen\yoff
\newcommand\clap[1]{\hbox to 0pt{\hss{#1}\hss}}
\newcommand\vdisk[1]{{\font\dotf=cmr10 scaled #1\dotf.}}
\newcommand\varline[2]{\setbox\dotb\hbox{\vdisk{#1}}\xoff=-.5\wd\dotb
\wd\dotb=0pt\yoff=-.5\ht\dotb\psep=#2\ht\dotb}
\newcommand\varpt[1]{\setbox\ptbox\clap{\vdisk{#1}}\setbox\ptbox
\hbox{\raise-.5\ht\ptbox\box\ptbox}}
\newcommand\cpt{\copy\ptbox}
\newcommand\point[3]{\rlap{\kern#1\unit\raise#2\unit\hbox{#3}}}
\newcommand\setnd[4]{\dx=#3\unit\advance\dx-#1\unit\divide\dx by\psep
\dy=#4\unit\advance\dy-#2\unit\divide\dy by\psep \multiply\dx
by\dx\multiply\dy by\dy\advance\dx\dy\nd=1\advance\dx-1sp
\loop\ifnum\dx>0\advance\dx-\nd sp\advance\nd1\advance\dx-\nd
sp\repeat}
\newcommand\dl[4]{{\setnd{#1}{#2}{#3}{#4}\dline{#1}{#2}{#3}{#4}\nd}}
\newcommand\dline[5]{{\nd=#5\hgt=#2\unit\dx=#3\unit\advance\dx-#1\unit
\divide\dx by\nd\dy=#4\unit\advance\dy-#2\unit\divide\dy by\nd
\advance\hgt\yoff\rlap{\kern#1\unit\kern\xoff\loop\ifnum\nd>1\advance\nd-1
\advance\hgt\dy\kern\dx\raise\hgt\copy\dotb\repeat}}}
\newcommand\ellipse[4]{\qellip{#1}{#2}{#3}{#4}\qellip{#1}{#2}{#3}{-#4}%
\qellip{#1}{#2}{-#3}{#4}\qellip{#1}{#2}{-#3}{-#4}}
\newcommand\qellip[4]{{\setnd{0}{0}{#3}{#4}\dx=\unit\dy=0pt\raise\yoff\rlap{%
\kern#1\unit\kern\xoff\raise#2\unit\hbox{\loop\ifnum\dx>0\rlap{\kern#3\dx
\raise#4\dy\copy\dotb}\hgt=\dx\divide\hgt
by\nd\advance\dy\hgt\hgt=\dy \divide\hgt
by\nd\advance\dx-\hgt\repeat\rlap{\raise#4\dy\copy\dotb}}}}}
\newcommand\bez[6]{{\setnd{#1}{#2}{#3}{#4}\ndx=\nd\setnd{#3}{#4}{#5}{#6}
\ifnum\ndx>\nd\nd=\ndx\fi\dx=#3\unit\advance\dx-#1\unit\dy=#4\unit
\advance\dy-#2\unit\dxx=#5\unit\advance\dxx-#1\unit\dyy=#6\unit\advance
\dyy-#2\unit\advance\dxx-2\dx\advance\dyy-2\dy\divide\dxx
by\nd\divide\dyy
by\nd\advance\dx.25\dxx\advance\dy.25\dyy\divide\dx
by\nd\divide\dy by\nd \multiply\nd
by2\dx=100\dx\dy=100\dy\dxx=100\dxx\dyy=100\dyy\divide\dxx by\nd
\divide\dyy
by\nd\hgt=#2\unit\raise\yoff\rlap{\kern#1\unit\kern\xoff
\raise\hgt\copy\dotb\loop\ifnum\nd>0\advance\nd-1\advance\hgt0.01\dy
\kern0.01\dx\raise\hgt\copy\dotb\advance\dx\dxx\advance\dy\dyy\repeat}}}
\newcommand\ptu[3]{\point{#1}{#2}{\cpt\raise1ex\clap{$\scriptstyle{#3}$}}}
\newcommand\ptd[3]{\point{#1}{#2}{\cpt\raise-1.8ex\clap{$\scriptstyle{#3}$}}}
\newcommand\ptr[3]{\point{#1}{#2}{\cpt\raise-.4ex\rlap{$\ \scriptstyle{#3}$}}}
\newcommand\ptl[3]{\point{#1}{#2}{\cpt\raise-.4ex\llap{$\scriptstyle{#3}\ $}}}
\newcommand\ptlu[3]{\point{#1}{#2}{\raise.8ex\clap{$\scriptstyle{#3}$}}}
\newcommand\ptld[3]{\point{#1}{#2}{\raise-1.6ex\clap{$\scriptstyle{#3}$}}}
\newcommand\ptlr[3]{\point{#1}{#2}{\raise-.4ex\rlap{$\,\scriptstyle{#3}$}}}
\newcommand\ptll[3]{\point{#1}{#2}{\raise-.4ex\llap{$\scriptstyle{#3}\,$}}}
\newcommand\pt[2]{\point{#1}{#2}{\cpt}}
\newcommand\thkline{\varline{1600}{.3}}
\newcommand\medline{\varline{800}{.5}}
\newcommand\thnline{\varline{400}{.6}}

\varpt{2500}\thnline\unit=2em

\newtheorem{thm}{Theorem}
                                                                                   
\newtheorem{conj}{Conjecture}
\newtheorem{qu}{Question}
\newtheorem{prob}{Problem}
\newtheorem{lemma}[thm]{Lemma}
\newtheorem{cor}[thm]{Corollary}

\newtheorem{obs}[thm]{Observation}
\theoremstyle{definition}\newtheorem{rmk}{Remark}
\theoremstyle{definition}
\newcommand{\ds}{\displaystyle}
\newcommand{\ul}{\underline}

\def\D{\mathcal{D}}

\def\F{\mathcal{F}}

\def\L{\mathcal{L}}

\def\P{\mathcal{P}}

\def\T{\mathcal{T}}

\def\X{\mathcal{X}}

\def\N{\mathbb{N}}

\def\RR{\mathbb{R}}

\def\le{\leqslant}
\def\ge{\geqslant}

\begin{document}
\title{Hereditary properties of tournaments}

\author{J\'ozsef Balogh}
\address{Department of Mathematics\\ University of Illinois\\ 1409 W. Green Street\\ Urbana, IL 61801} \email{jobal@math.uiuc.edu}

\author{B\'ela Bollob\'as}
\address{Department of Mathematical Sciences\\ The University of Memphis\\ Memphis, TN 38152\\ and\\ Trinity College\\ Cambridge CB2 1TQ\\ England} \email{B.Bollobas@dpmms.cam.ac.uk}

\author{Robert Morris}
\address{Department of Mathematical Sciences\\ The University of Memphis\\ Memphis, TN 38152} \email{rdmorrs1@memphis.edu}\thanks{The first author was supported during this research by OTKA grant T049398 and NSF grant DMS-0302804, the second by NSF grant CCR-0225610 and ARO grant W911NF-06-1-0076, and the third by a Van Vleet Memorial Doctoral Fellowship.}

\begin{abstract}
A collection of unlabelled tournaments $\P$ is called a \emph{hereditary property} if it is closed under isomorphism and under taking induced sub-tournaments. The speed of $\P$ is the function $n \mapsto |\P_n|$, where $\P_n = \{T \in \P : |V(T)| = n\}$. In this paper, we prove that there is a jump in the possible speeds of a hereditary property of tournaments, from polynomial to exponential speed. Moreover, we determine the minimal exponential speed, $|\P_n| = c^{(1+o(1))n}$, where $c \simeq 1.47$ is the largest real root of the polynomial $x^3 = x^2 + 1$, and the unique hereditary property with this speed.
\end{abstract}

\maketitle

\section{Introduction}\label{tournintro}

In this paper we shall prove that there is a jump in the possible speeds of a hereditary property of tournaments, from polynomial to exponential speed. We shall also determine the minimum possible exponential speed, and the unique hereditary property giving rise to this speed. This minimum speed is different from those previously determined for other structures (see \cite{order}, \cite{KK}, \cite{Klaz}). In order to state our result, we shall need to begin with some definitions.

A tournament is a complete graph with an orientation on each edge. Here we shall deal with \emph{unlabelled} tournaments, so two tournaments $S$ and $T$ are isomorphic if there exists a bijection $\phi: V(S) \to V(T)$ such that $u \to v$ if and only if $\phi(u) \to \phi(v)$. Throughout the paper, we shall not distinguish isomorphic tournaments. A property of tournaments is a collection of unlabelled tournaments closed under isomorphisms of the vertex set, and a property of tournaments is called \emph{hereditary} if it is closed under taking sub-tournaments. If $\P$ is a property of tournaments, then $\P_n$ denotes the collection $\{T \in \P : |V(T)| = n\}$, and the function $n \mapsto |\P_n|$ is called the \emph{speed} of $\P$. Analogous definitions can be made for other combinatorial structures (e.g., graphs, ordered graphs, posets, permutations).

We are interested in the (surprising) phenomenon, observed for hereditary properties of various types of structure (see for example~\cite{AlekS}, \cite{BT}, \cite{MT}) that the speeds of such a property are far from arbitrary. More precisely, there often exists a family $\F$ of functions $f : \N \to \N$ and another function $F : \N \to \N$, with $F(n)$ {\em much} larger than $f(n)$ for every $f \in \F$, such that if for each $f \in \F$ the speed is infinitely often larger than $f(n)$, then it is also larger than $F(n)$ for every $n \in \N$. Putting it concisely: the speed {\em jumps} from $\F$ to $F$.

Hereditary properties of \emph{labelled} oriented graphs, and in particular properties of posets, have been extensively studied. For example, Alekseev and Sorochan~\cite{AlekS} proved that the labelled speed $|\P^n|$ of a hereditary property of oriented graphs is either $2^{o(n^2)}$, or at least $2^{n^2/4 + o(n^2)}$, and Brightwell, Grable and Pr\"omel~\cite{BGP} showed that for a principal hereditary property of labelled posets (a property in which only one poset is forbidden), either $(i) \; |\P^n| \le n! \, c^n$ for some $c \in \RR$, $(ii) \; n^{c_1n} \le |\P^n| \le n^{c_2n}$ for some $c_1,c_2 \in \RR$, $(iii) \; n^{Cn} \ll |\P^n| = 2^{o(n^2)}$ for every $C \in \RR$, or $(iv) \; |\P^n| = 2^{n^2/4 + o(n^2)}$. Other papers on the speeds of particular poset properties include~\cite{AS}, \cite{BG} and \cite{BPS}. For properties of labelled graphs, Balogh, Bollob\'as and Weinreich~\cite{BBW1}, \cite{BBW2} have determined the possible speeds below $n^{n + o(n)}$ very precisely, and their proofs can be adapted to prove corresponding results for labelled oriented graphs. Much is still unknown, however, about properties with speed $|\P^n| = n^{n + o(n)}$, and about those with speed greater than $2^{n^2/4}$.

For very high speed properties, the unlabelled case is essentially the same as the labelled case, since the speeds differ by a factor of only at most $n!$ (the total possible number of labellings). However, for properties with lower speed, the two cases become very different, and the unlabelled case becomes much more complicated. For example, the speed of a hereditary property of labelled tournaments is either zero for sufficiently large $n$, or at least $n!$ for every $n \in \N$, simply because any sufficiently large tournament contains a transitive sub-tournament on $n$ vertices (see Observation~\ref{trans}), and such a tournament is counted $n!$ times in $|\P^n|$. On the other hand, in \cite{posets} the authors found it necessary to give a somewhat lengthy proof of the following much smaller jump: a hereditary property of unlabelled tournaments has either bounded speed, or has speed at least $n - 2$.

Given the difficulty we had in proving even this initial jump, one might suspect that describing all polynomial-speed hereditary properties of tournaments, or proving a jump from polynomial to exponential speed for such properties, would be a hopeless task. However, in Theorem~\ref{tourns} (below) we shall show that this is not the case. Indeed, we shall prove that the speed $|\P_n|$ of a hereditary property of unlabelled tournaments is either bounded above by a polynomial, or is at least $c^{(1+o(1))n}$, where $c \simeq 1.47$ is the largest real root of the polynomial $x^3 = x^2 + 1$.

Results analogous to Theorem~\ref{tourns} have previously been proved for labelled graphs~\cite{BBW1}, labelled posets~\cite{posets}, permutations~\cite{KK} and ordered graphs~\cite{order}. In the latter two cases the minimum exponential speed is the sequence $F_n = F_{n-1} + F_{n-2}$, the Fibonacci numbers, and in \cite{Klaz} there was an attempt to characterize the structures whose growth admits Fibonacci-type jumps. Tournament properties were not included in this characterization and, as Theorem~\ref{tourns} shows, they exhibit a different (though similar) jump from polynomial to exponential speed.

Each of these results is heavily dependent on the labelling/order on the vertices. When dealing with unlabelled and unordered vertices, we have many possible isomorphisms to worry about (instead of only one), so many new problems are created. The only result similar to Theorem~\ref{tourns} for such structures, of which we are aware, is for unlabelled graphs~\cite{BBSS}. The proof in that paper uses the detailed structural results about properties of labelled graphs proved in~\cite{BBW1} and \cite{BBW2}; in contrast, our proof is self-contained.

\section{Main Results}

In this section we shall state our main results. We begin by describing the hereditary property of tournaments with minimal exponential speed.

Consider the following collection $\T$ of tournaments. For each $m \in \N$ and $a_1, \ldots, a_m \in \{1,3\}$, let $T = T(a_1, \ldots, a_m)$ be the tournament with vertex set $\{x(i,j) : i \in [m], j \in [a_i]\}$, in which
$$x(i,j) \to x(k,\ell)\textup{ if }i < k\textup{, or if }i = k\textup{  and }k  - i \equiv 1 \hspace{-0.25cm} \pmod 3.$$
Thus $|V(T)| = \sum_{i=1}^m a_i$, and the sequence $(a_1,\ldots,a_m)$ can be reconstructed from $T$ (see Lemma~\ref{I=111}). Define $$\T \; = \; \{T(a_1, \ldots, a_m) \,:\, m \in \N, \: a_1, \ldots, a_m \in \{1,3\}\},$$ and note that $\T$ is a hereditary property of tournaments.

Now, let $F^*_n$ be the Fibonacci-type sequence of integers defined by $F^*_0 = F^*_1 = F^*_2 = 1$, and $F^*_n = F^*_{n-1} + F^*_{n - 3}$ for every $n \ge 3$. Note that $F^*_n = c^{(1+o(1))n}$ as $n \to \infty$, where $c \simeq 1.47$ is the largest real root of the polynomial $x^3 = x^2 + 1$. Note also that $|\T_n| = F_n^*$ for every $n \in \N$ (again, see Lemma~\ref{I=111} for the details).

The following theorem, which is the main result of this paper, says that $\T$ is the unique smallest hereditary property of tournaments with super-polynomial speed.

\begin{thm}\label{tourns}
Let $\P$ be a hereditary property of tournaments. Then either
\begin{enumerate}
\item[$(a)$] $|\P_n| = \Theta(n^k)$ for some $k \in \N$, or\\[-1ex]
\item[$(b)$] $|\P_n| \ge F^*_n$ for every $4 \neq n \in \N$.\\[-2.5ex]
\end{enumerate}
Moreover, this lower bound is best possible, and $\T$ is the unique hereditary property of tournaments $\P$ with $|\P_n| = F^*_n$ for every $n \in \N$.
\end{thm}

We remark that there exists a hereditary property of tournaments $\P$ with speed roughly $2^n$, but for which $|\P_4| = 2 < 3 = F_4^*$ (see Lemma~\ref{osmall}), so this result really is best possible.

Our second theorem determines the speed of a polynomial-speed hereditary property of tournaments up to a constant. The statement requires the notion of a \emph{homogeneous block} in a tournament, which will be defined in Section~\ref{tkeysec}, but we state it here in any case, for ease of reference. Given a hereditary property of tournaments $\P$, let $k(\P) = \sup\{ \ell : \forall \, m \in \N, \: \exists \, T \in \P$ such that the $(\ell + 1)^{st}$ largest homogeneous block in $T$ has at least $m$ elements$\}$.

\begin{thm}\label{tpolythm}
Let $\P$ be a hereditary property of tournaments, and suppose that $k = k(\P) < \infty$. Then $|\P_n| = \Theta(n^k)$.
\end{thm}

The proof of Theorem~\ref{tourns} is roughly as follows. In Section~\ref{tkeysec} we shall define the homogeneous block decomposition of a tournament, and show that if the number of distinct homogeneous blocks occurring in a tournament in $\P$ is bounded, then the speed of $\P$ is bounded above by a polynomial, whereas if it this number is unbounded, then certain structures must occur in $\P$. Then, in Section~\ref{T1T2}, we shall use the techniques developed in \cite{order} to show that if these structures occur, then the speed must be at least $F^*_n$. In Section~\ref{tpolysec} we shall investigate the possible polynomial speeds, and prove Theorem~\ref{tpolythm}, and in Section~\ref{tproof} we put the pieces together and prove Theorem~\ref{tourns}. In Section~\ref{tprobs} we shall discuss possible future work.

We shall use the following notation throughout the paper. If $n \in \N$ and $A, B \subset \N$, we say that $n > A$ if $n > a$ for every $a \in A$, and $A > B$ if $a > b$ for every $a \in A$ and $b \in B$. Also, if $T$ is a tournament, $v \in V(T)$ and $C, D \subset V(T)$, then we say that $v \to C$ if $v \to c$ for every $c \in C$, and $C \to D$ if $c \to d$ for every $c \in C$, $d \in D$. We shall sometimes write $u \in T$ to mean that $u$ is a vertex of $T$. Finally, $[n] = \{1, \ldots, n\}$, and $[0] = \emptyset$.

\section{Homogeneous blocks, and the key lemma}\label{tkeysec}

We begin by defining the concept of a homogeneous block in a tournament. Let $T$ be a tournament, and let $u,v \in T$. Write $u \curvearrowright v$ if $u \to v$, and for some $k \ge 0$ and some set of vertices $w_1, \ldots, w_k$, the following conditions hold. Let $C(u,v) = \{u, w_1, \ldots, w_k, v\}$.\\[-2ex]
\begin{enumerate}
\item[$(i)$] $u \to w_i \to v$ for every $1 \le i \le k$,\\[-1.5ex]
\item[$(ii)$] $w_i \to w_j$ for every $1 \le i < j \le k$, and\\[-1.5ex]
\item[$(iii)$] if $x \in V(T) \setminus C(u,v)$, and $y,z \in C(u,v)$, then $x \to y$ if and only if $x \to z$.\\[-2ex]
\end{enumerate}
We say that the pair $\{u,v\}$ is \emph{homogeneous} (and write $u \sim v$) if $u = v$, or $u \curvearrowright v$, or $v \curvearrowright u$. If $u \curvearrowright v$, then we call $C(u,v)$ the \emph{homogeneous path} from $u$ to $v$, and define $C(v,u) = C(u,v)$. Note that $C(u,v)$ is well-defined, since if it exists (and $u \to v$, say), then it is the set $\{u,v\} \cup \{w : u \to w \to v\}$. Note also that $x \sim y$ for every pair $x,y \in C(u,v)$.

\begin{lemma}\label{equiv}
$\sim$ is an equivalence relation.
\end{lemma}

\begin{proof}
Symmetry and reflexivity are clear; to show transitivity, consider vertices $x$, $y$ and $z$ in $T$ with $x \sim y$ and $y \sim z$, and suppose without loss that $x \to y$. We shall show that $x \sim z$. Let the sets $C(x,y) = \{x, w_1, \ldots, w_k, y \}$ and $C(y,z) = \{y, w'_1, \ldots, w'_\ell, z \}$ be the homogeneous paths from $x$ to $y$ and between $y$ and $z$ respectively. As noted above, if $z \in C(x,y)$ then $x \sim z$, so we are done, and similarly if $x \in C(y,z)$ then $x \sim z$.

So assume that $z \notin C(x,y)$ and $x \notin C(y,z)$. Now if $z \to y$, then also $z \to x$, since $x \sim y$ and $z \notin C(x,y)$. But then $z \to x \to y$, so $x \in C(y,z)$, a contradiction. Hence $y \to z$, and so $C(x,y) \cap C(y,z) = \{y\}$, since $w \to y$ if $y \neq w \in C(x,y)$ and $y \to w$ if $y \neq w \in C(y,z)$. But now $x \sim z$, with $C(x,z) = C(x,y) \cup C(y,z)$, since for any $w \in C(x,y)$ and $w' \in C(y,z)$ with $w \neq w'$, $y \to w'$ so $w \to w'$, and for any $v \notin C(x,y) \cup C(y,z)$ and any $w,w' \in C(x,y) \cup C(y,z)$, $w \to v$ if and only if $y \to v$, if and only if $w' \to v$.
\end{proof}

We may now define a \emph{homogeneous block} in a tournament $T$ to be an equivalence class of the relation $\sim$. By Lemma~\ref{equiv}, we may partition the vertices of any tournament $T$ into \emph{homogeneous blocks} in a unique way (see Figure 1).\\

\[ \unit = 0.5cm \hskip -13\unit
\medline \ellipse{6}{0}{3}{1} \ellipse{0}{-5}{1}{3} \ellipse{6}{-10}{3}{1} \ellipse{12}{-5}{1}{3}
\point{0}{-3.5}{$ \dl{0}{0}{0}{-3} \dl{0}{-3}{0.3}{-2.5} \dl{0}{-3}{-0.3}{-2.5} $}
\point{12}{-3.5}{$ \dl{0}{0}{0}{-3} \dl{0}{0}{0.3}{-0.5} \dl{0}{0}{-0.3}{-0.5} $}
\point{0}{0}{$ \dl{4.5}{0}{7.5}{0} \dl{4.5}{0}{5}{0.3} \dl{4.5}{0}{5}{-0.3} $}
\point{0}{-10}{$ \dl{4.5}{0}{7.5}{0} \dl{7.5}{0}{7}{0.3} \dl{7.5}{0}{7}{-0.3} $}
\dl{2.7}{-0.7}{1}{-2} \dl{1}{-2}{1.6}{-1.95} \dl{1}{-2}{1.205}{-1.435}
\dl{9.3}{-9.3}{11}{-8} \dl{11}{-8}{10.4}{-8.05} \dl{11}{-8}{10.795}{-8.565}
\dl{11}{-2}{9.3}{-0.7} \dl{9.3}{-0.7}{9.9}{-0.75} \dl{9.3}{-0.7}{9.505}{-1.265}
\dl{1}{-8}{2.7}{-9.3} \dl{2.7}{-9.3}{2.1}{-9.25} \dl{2.7}{-9.3}{2.495}{-8.735}
\bez{2.5}{-5}{6}{-3}{9.5}{-5} \dl{9.5}{-5}{8.9}{-5} \dl{9.5}{-5}{9.145}{-4.485}
\dl{6}{-7.5}{6}{-2.5} \dl{6}{-7.5}{6.3}{-7} \dl{6}{-7.5}{5.7}{-7}
\point{0.5}{-13}{Figure 1: Homogeneous blocks}
\]\\

Let $B(T)$ denote the number of homogeneous blocks of a tournament $T$, and if $\P$ is a property of tournaments, let $B(\P)$ denote $\sup\{B(T) : T \in \P\}$, where $B(\P)$ may of course be equal to infinity.

\begin{lemma}\label{Bbdd}
Let $\P$ be a hereditary property of tournaments, and let $M \in \N$. If $B(\P) = M + 1$, then $|\P_n| = O(n^M)$.
\end{lemma}

\begin{proof}
Let $B(\P) = M+1$. Each tournament $T \in \P_n$ is determined by a sequence $(a_1, \dots, a_{M+1})$ of non-negative integers summing to $n$, and an ordered tournament on $M+1$ vertices. Thus $$\ds|\P_n| \le 2^{(M+1)^2}{{n+M} \choose M} = O(n^M),$$ as claimed.
\end{proof}

We shall now prove the key lemma in the proof of Theorem~\ref{tourns}. We first need to define some particular structures, which will play a pivotal role in the proof; they come in two flavours. Let $T$ be a tournament, and let $k \in \N$.
\begin{itemize}
\item Type~1: there exist distinct vertices $x_1, \ldots, x_{2k}$ and $y$ in $T$ such that $x_i \to x_j$ if $i < j$, and $y \to x_i$ if and only if $x_{i+1} \to y$, for each $i \in [2k-1]$.\\[-1.5ex]
\item Type $2$: there exist distinct vertices $x_1, \ldots, x_{2k}$ and $y_1, \ldots, y_k$ in $T$ such that $x_i \to x_j$ if $i < j$, and $x_{2i} \to y_i \to x_{2i-1}$ for every $i \in [k]$.
\end{itemize}

Note that there are two different structures of Type~1, and only one of Type 2. We refer to these as $k$-structures of Type~1 and 2. Type 2 structures are not tournaments, but sub-structures contained in tournaments: instead of saying that ``a structure of Type $2$ occurs in $\P$'' it would be more precise to say that ``there is a tournament $T \in \P$ admitting a structure of Type $2$". However, for smoothness of presentation we sometimes handle them as tournaments.

We shall use the following simple observation, which may easily proved by induction.

\begin{obs}\label{trans}
A tournament on at least $2^n$ vertices contains a transitive subtournament on at least $n$ vertices.
\end{obs}

The following lemma is the key step in the proof of Theorem~\ref{tourns}.

\begin{lemma}\label{tournkey}
Let $\P$ be a hereditary property of tournaments. If $B(\P) = \infty$, then $\P$ contains arbitrarily large structures of Type~$1$ or~$2$.
\end{lemma}

\begin{proof}
Let $\P$ be a hereditary property of tournaments with $B(\P) = \infty$, and let $k \in \N$. We shall show that $\P$ contains either a $k$-structure of Type~1, or a $k$-structure of Type 2 (or both).

To do this, first let $K = 4k^2 2^{16k^6} + 8k^2$, let $M = 2^K$, and let $T' \in \P$ be a tournament with at least $M$ different homogeneous blocks. Choose one vertex from each block, and let $T$ be the tournament induced by those vertices. Note that $T \in \P$, and that the homogeneous blocks of $T$ are single vertices, since if $x \sim y$ in $T$, then $x \sim y$ in $T'$. Thus, for each pair of vertices $x,y \in V = V(T)$, there exists a vertex $z \in V$ such that $x \to z \to y$ or $y \to z \to x$.

Let $A$ be the vertex set of a maximal transitive sub-tournament of $T$, so by Observation~\ref{trans}, $|A| = r \ge K$. Order the vertices of $A = \{a_1, \ldots, a_r\}$ so that $a_i \to a_j$ if $i < j$. Then, for each pair $\{a_i,a_{i+1}\}$ with $i \in [r - 1]$, choose a vertex $b_i \in V \setminus A$ such that $a_{i+1} \to b_i \to a_i$ if one exists; otherwise choose $b_i$ such that $a_i \to b_i \to a_{i+1}$. As observed above, such a $b_i$ must exist. Let $Y = \{b_i : a_{i+1} \to b_i \to a_i\}$, and for each $y \in Y$, let $Z_y = \{a_i \in A : a_{i+1} \to y \to a_i\}$. The following two claims show that $a_{i+1} \to b_i \to a_i$ for only a bounded number of indices $i$.\\

\noindent\ul{Claim 1}: If $|Z_y| \ge 2k$ for some $y \in Y$, then $T$ contains a $k$-structure of Type~1.

\begin{proof}
Let $y \in Y$ and suppose that $|Z_y| \ge 2k$. Let $Z' = \{a_{i(1)}, \ldots, a_{i(2k)}\}$ be any subset of $Z_y$ of order $2k$, and suppose $i(1) < \ldots < i(2k)$. By definition, $a_{i(j) + 1} \to y \to a_{i(j)}$ for each $j \in [2k]$. Let $Z'' = \{a_{i(2j - 1)}, a_{i(2j - 1) + 1} : j \in [k]\}$. Since $T[A]$ is transitive, so is $T[Z'']$, and thus $T[Z'' \cup \{y\}]$ is a $k$-structure of Type~1.
\end{proof}

\noindent\ul{Claim 2}: If $|Y| \ge 2k$ then $T$ contains a $k$-structure of Type 2.

\begin{proof}
Suppose $|Y| \ge 2k$, and let $Y'$ be any subset of $Y$ of order $2k$. For each vertex $y \in Y'$, choose an index $i = i(y)$ such that $y = b_i$, and note that $i(y) = i(y')$ implies $y = y'$. Let $I = \{i(y) : y \in Y'\}$ have elements $i_1 < \ldots < i_{2k}$, and let $I' = \{i_1,i_3,\ldots, i_{2k-1}\}$. Finally, let $A' = \{a_i \in A : i \in I'$ or $i-1 \in I'\}$, and let $Y'' = \{y \in Y' : y = b_i$ for some $i \in I'\}$, so $i(y) \in I'$ if $y \in Y''$.

We claim that $T[A' \cup Y'']$ contains a $k$-structure of Type 2. Indeed, $T[A']$ is transitive (since $A$ is transitive), and if $y \in Y''$, then $a_{i(y)}, a_{i(y)+1} \in A'$, and $a_{i(y)+1} \to y \to a_{i(y)}$ in $T$. Moreover, since we used only every other entry of $I$, no two of the pairs $\{a_{i(y)},a_{i(y)+1}\}$ overlap. Thus $|A'| = 2k$ and $|Y''| = k$, so $T[A' \cup Y'']$ contains a $k$-structure of Type 2, as claimed.
\end{proof}

If $|Y| \ge 2k$, or if $|Z_y| \ge 2k$ for any $y \in Y$, then we are done by Claims 1 and 2. So assume that $|Y| < 2k$ and that $|Z_y| < 2k$ for every $y \in Y$. Let $$P = \{\{a_i,a_{i+1}\} \subset A : a_{i+1} \to v \to a_i\textup{ for some }v \in V \setminus A\}$$ be the set of consecutive pairs of $A$ which are contained in some cyclic triangle of $T$. Since we chose $b_i$ such that $a_{i+1} \to b_i \to a_i$ if possible, each pair in $P$ contributes one vertex to $Z_y$ for at least one $y \in Y$. Thus $|P| \le \ds\sum_{y \in |Y|} |Z_y| < 4k^2$. Therefore, by the pigeonhole principle, there must exist an interval $C \subset [r - 1]$ of size at least $(r - 8k^2)/4k^2 \ge 2^{16k^6}$, such that $A_C = \{a_i : i \in C\}$ contains no element of any pair of $P$. In other words, $\{i,i+1\} \cap C = \emptyset$ for every pair $\{a_i,a_{i+1}\} \in P$.

Now, for each $i \in C$, recall that $b_i \in V \setminus A$, the vertex chosen earlier, satisfies $a_i \to b_i \to a_{i+1}$. Let $X = \{b_j : j \in C\}$. Observe that $a_i \to b_j \to a_{i'}$ for every $i,j,i' \in C$ with $i \le j < i'$, since otherwise there must exist a pair of consecutive vertices $a_\ell$ and $a_{\ell+1}$ of $A$, with $\ell \in C$, such that $a_{\ell+1} \to b_j \to a_\ell$, contradicting the definition of $C$. Hence the vertices $b_j$ with $j \in C$ are all distinct.

It follows that $|X| = |C| \ge 2^{16k^6}$. Therefore, by Observation~\ref{trans}, there exists a transitive sub-tournament of $T[X]$ on $s \ge 16k^6$ vertices. Let the vertex set of this transitive sub-tournament be $X' = \{x(1), \ldots, x(s)\}$, ordered so that $x(i) \to x(j)$ if $i < j$, and let $C' = \{i \in C : b_i \in X'\}$.

Define $\phi: C' \to [s]$ to be the function such that $x(\phi(i)) = b_i$. Note that $\phi$ is surjective. By the Erd\H{o}s-Szekeres Theorem, there exists a subset $C''$ of $C'$ of order $t \ge \sqrt{s} \ge 4k^3$, such that $\phi$ is either strictly increasing or strictly decreasing on $C''$. Let $X'' = \{ b_i \in X' : i \in C''\}$ be the corresponding subset of $X'$.

The following two claims now complete the proof of the lemma.\\

\noindent\ul{Claim 3}: If $\phi$ is increasing on $C''$, then $T$ contains a $k$-structure of Type~1.

\begin{proof}
Suppose that $\phi$ is strictly increasing on $C''$, so $b_i \to b_j$ for every $i,j \in C''$ with $i < j$. Recall also that $a_i \to b_j \to a_{i'}$ for every $i,j,i' \in C$ with $i \le j < i'$. Thus $T[\{a_i, b_i : i \in C''\}]$ is a transitive tournament, with $a_i \to b_i \to a_j \to b_j$ for every $i,j \in C''$ with $i < j$.

Let $v \in V \setminus A$. Since $A$ is a maximal transitive sub-tournament, $T[A \cup \{v\}]$ is not transitive, so $a_{i+1} \to v \to a_i$ for some $i \in [r-1]$. Recall that $|P| < 4k^2$, so by the pigeonhole principle, there must exist a consecutive pair $\{a_\ell,a_{\ell+1}\} \in P$, and a subset $W \subset X''$ of order $q \ge |X''|/4k^2 \ge k$, such that $a_{\ell+1} \to w \to a_\ell$ for each vertex $w \in W$. Note that (by the definition of $C$) either $\ell + 1 < C$ or $\ell > C$.

Now, let $D = \{i \in C : b_i \in W\}$ be the subset of $C''$ corresponding to $W$, with elements $d(1) < \ldots < d(q)$. Let $E = \{a_i : i \in D\}$ be the corresponding subset of $A$. Then, by the comments above, $T[E \cup W]$ is a transitive tournament with vertices $a_{d(1)} \to b_{d(1)} \to \ldots \to a_{d(q)} \to b_{d(q)}$.

Suppose $\ell + 1 < C$, where $\{a_\ell, a_{\ell+1}\}$ is the pair defined earlier. Then $a_\ell \to a_{d(i)}$ and $b_{d(i)} \to a_\ell$ for every $i \in [q]$. It follows that $T[\{a_\ell\} \cup E \cup W]$ is a $q$-structure of Type~1, and so contains a $k$-structure of Type~1 (since $q \ge k$). Similarly, if $\ell > C$, then $a_{d(i)} \to a_{\ell+1}$ and $a_{\ell+1} \to b_{d(i)}$ for every $i \in [q]$. It follows that $T[\{a_{\ell+1}\} \cup E \cup W]$ is a $q$-structure of Type~1, and again we are done.
\end{proof}

\noindent\ul{Claim 4}: If $\phi$ is decreasing on $C''$, then $T$ contains a $k$-structure of Type 2.

\begin{proof}
Suppose that $\phi$ is strictly decreasing on $C''$, so $b_j \to b_i$ for every $i,j \in C''$ with $i < j$. Let $C'' = \{c(1), \ldots, c(t)\}$, with $c(1) < \ldots < c(t)$. As in Claim 3, we have $a_{c(i)} \to b_{c(i)} \to a_{c(i+1)}$ for every $i \in [t-1]$.

Now $T[X'']$ is transitive, with $b_{c(t)} \to \ldots \to b_{c(1)}$. Also $b_{c(2i-1)} \to a_{c(2i)} \to b_{c(2i)}$ for every $i \in [t^*]$, where $t^* = \left\lfloor t/2 \right\rfloor$. Thus, letting $C''' = \{c(2i) : i \in [t^*]\}$, we have shown that $T[X'' \cup C''']$ contains a $t^*$-structure of Type 2. Since $t \ge 4k^3 \ge 2k$, this proves the claim.
\end{proof}

By the comments above, $\phi$ is either strictly increasing or strictly decreasing, so this completes the proof of the lemma.
\end{proof}

Combining Lemmas~\ref{Bbdd} and \ref{tournkey}, we get the following result, which summarises what we have proved so far.

\begin{cor}\label{tsumm}
Let $\P$ be a hereditary property of tournaments. If for every $k \in \N$ there are  infinitely many values of $n$ such that $|\P_n| \ge n^k$, then $\P$ contains arbitrarily large structures of Type~1 or 2.
\end{cor}

\section{Structures of Type~1 and 2}\label{T1T2}

We begin by showing that if $\P$ contains arbitrarily large structures of Type~1, then the speed of $\P$ is at least $2^{n-1} - O(n^2)$. Given $n \in \N$, and a subset $S \subset [n]$, let $T_{n+1}(S)$ denote the tournament with $n + 1$ vertices, $\{y, x_1, \ldots, x_n\}$ say, in which $x_i \to x_j$ if $i < j$, and $y \to x_i$ if and only if $i \in S$. Suppose $T = T_{n+1}(S)$ has exactly one transitive sub-tournament on $n$ vertices (i.e., there is exactly one subset $A \subset V(T)$ with $|A| = n$ such that $T[A]$ is transitive). Then the vertex $y \in V(T)$ is uniquely determined, and so the set $S$ is determined by $T$. Hence if $S$ and $S'$ are distinct sets satisfying that $T_{n+1}(S)$ and $T_{n+1}(S')$ each have exactly one transitive sub-tournament on $n$ vertices, then $T_{n+1}(S)$ and $T_{n+1}(S')$ are distinct tournaments.

For each $n \in \N$, let $\D_n = \{S \subset [n] : T_{n+1}(S)$ has at least two transitive sub-tournaments on $n$ vertices$\}$. We shall use the following simple observation to prove Lemma~\ref{T1}.

\begin{obs}\label{sizeD}
$|\D_n| \le 2\ds{n \choose 2} + n + 1$ for every $n \in \N$.
\end{obs}

\begin{proof}
Let $n \in \N$, $S \subset [n]$, and suppose that $S \in \D_n$. Let $T = T_{n+1}(S)$ have vertex set $V = \{y, x_1, \ldots, x_n\}$ say, where $x_i \to x_j$ if $i < j$, and $y \to x_i$ if and only if $i \in S$. Let $\ell \in [n]$ be such that $T - \{x_\ell\}$ is transitive (such an $\ell$ exists because $S \in \D_n$).

Now, $T - \{x_\ell\}$ is transitive, so there exists an $m \in [0,n]$ such that $x_i \to y$ if $\ell \neq i \le m$, and $y \to x_i$ if $\ell \neq i > m$. There are three cases to consider.\\[-1.5ex]

\noindent \textbf{Case 1}: If $m \le \ell$ and $y \to x_\ell$, or $m \ge \ell - 1$ and $x_\ell \to y$, then $T$ is transitive, and $S = [i,n]$ for some $i \in [n+1]$.\\[-1.5ex]

\noindent \textbf{Case 2}: If $m \ge \ell + 1$ and $y \to x_\ell$, then $S = \{\ell\} \cup [m+1,n]$.\\[-1.5ex]

\noindent \textbf{Case 3}: If $m \le \ell - 2$ and $x_\ell \to y$, then $S = [m+1,\ell-1] \cup [\ell+1,n]$.\\[-1.5ex]

So the set $S$ must be of the form $[i,n]$ with $i \in [n+1]$, or $\{i\} \cup [j+1,n]$ with $1 \le i < j \le n$, or $[i,j-1] \cup [j+1,n]$ with $1 \le i < j \le n$. There are at most $\ds2{n \choose 2} + n + 1$ such sets.
\end{proof}

The following lemma gives the desired lower bound when $\P$ contains arbitrarily large structures of Type~1.

\begin{lemma}\label{T1}
Let $\P$ be a hereditary property of tournaments. Suppose $k$-structures of Type~1 occur in $\P$ for arbitrarily large values of $k$. Then
$$|\P_n| \; \ge \; 2^{n-1} - 2\ds{{n-1} \choose 2} - n$$ for every $n \in \N$.
\end{lemma}

\begin{proof}
Let $\P$ be a hereditary property of tournaments containing arbitrarily large structures of Type~1, and let $n \in \N$. Let $T \in \P$ be an $n$-structure of Type~1, with vertex set $V = \{y,x_1, \ldots, x_{2n}\}$, where $x_i \to x_j$ if $i < j$, and $y \to x_i$ if and only if $i$ is odd. For $n = 1$ the result is trivial, so assume that $n \ge 2$.

We claim that $T_n(S)$ is a sub-tournament of $T$ for every $S \subset [n-1]$. Indeed, let $S \subset [n-1]$, and let $T'$ be induced by the vertices $y \cup \{x_{2i-1} : i \in S\} \cup \{x_{2i} : i \notin S\}$. It is easy to see that $T' = T_n(S)$.

Now, every sub-tournament of $T$ has some transitive sub-tournament on $n-1$ vertices. As noted above, if $T_n(S)$ and $T_n(S')$ each have exactly one transitive sub-tournament on $n-1$ vertices, and $S \neq S'$, then they are distinct. Hence the tournaments $\{T_n(S) : S \subset [n-1], S \notin \D_{n-1}\} \subset \P_n$ are all distinct. By Observation~\ref{sizeD}, there are at least $2^{n-1} - \ds{{n-1} \choose 2} - n$ subsets $S \subset [n-1]$, $S \notin \D_{n-1}$. The result follows.
\end{proof}

\begin{rmk}
With a little more work one can replace the lower bound in Lemma~\ref{T1} by $|\P_n| \ge 2^{n-1} - \ds{{n-1} \choose 2} - 1$, which is best possible for $n \ge 2$. We shall not need this sharp result however.
\end{rmk}

We now turn to those properties containing arbitrarily large structures of Type 2. We shall show that such a property has speed at least $F^*_n$. We first need to define eight specific (families of) tournaments. We shall show that if $\P$ contains arbitrarily large structures of Type 2, then it contains arbitrarily large members of one of these families.

Given $n \in \N$, let $X = \{x_1, \ldots, x_{2n}\}$ and $Y = \{y_1, \ldots, y_n\}$ be disjoint ordered sets of vertices, and let $I = (I_1, I_2, I_3) \in \{0,1\}^3$. Define $M_I(X,Y) = M_I^{(n)}$ to be the tournament with vertex set $X \cup Y$, and with edges oriented as follows.
\begin{enumerate}
\item[$(i)$] \;If $1 \le i < j \le 2n$ then $x_i \to x_j$,
\item[$(ii)$] \;for each $1 \le i \le n$, $x_{2i} \to y_i \to x_{2i-1}$,
\item[$(iii)$] \;if $1 \le i < j \le n$, then $y_i \to y_j \Leftrightarrow I_1 = 1$,
\item[$(iv)$] \;if $i \in [2n]$, $j \in [n]$ and $i \le 2j - 2$, then $x_i \to y_j \Leftrightarrow I_2 = 1$,
\item[$(v)$] \;if $i \in [2n]$, $j \in [n]$ and $i \ge 2j + 1$, then $y_j \to x_i \Leftrightarrow I_3 = 1$.
\end{enumerate}
We remark that the family $\{T : T \le M_{(1,1,1)}^{(n)}$ for some $n \in \N\}$ is exactly the family $\T$ defined in the introduction.

\begin{lemma}\label{t2}
Let $\P$ be a hereditary property of tournaments, and let $K \in \N$. Suppose that no $k$-structures of Type~1 occur in $\P$ for $k \ge K$, but that $k$-structures of Type 2 occur in $\P$ for arbitrarily large values of $k$. Then $M_I^{(n)} \in \P$ for some $I \in \{0,1\}^3$, and every $n \in \N$.
\end{lemma}

\begin{proof}
We shall use Ramsey's Theorem. Recall that $R_r(s)$ denotes the smallest number $m$ such that any $r$-colouring of the edges of $K_m$ contains a monochromatic $K_s$. Let $K \in \N$, and $\P$ be a hereditary property of tournaments as described, and let $n \in \N$, $R = R_{16}(\max\{n, K+1\})$, and $k = 2^{R^2}$. Choose a tournament $T \in \P$ containing a $k$-structure of Type 2 on vertices $\{x_1, \ldots, x_{2k}, y_1, \ldots, y_k\}$, so $x_i \to x_j$ if $i < j$, and $x_{2i} \to y_i \to x_{2i-1}$ for every $i \in [k]$.

First, by Observation~\ref{trans}, there exists a subset $Y \subset \{y_1, \ldots, y_k\}$, with $|Y| \ge R^2$, such that $T[Y]$ is transitive. Now, by the Erd\H{o}s-Szekeres Theorem, there exists a subset $Y' \subset Y$, with $|Y'| = t \ge \sqrt{|Y|} \ge R$, such that either $y_i \to y_j$ for every $y_i, y_j \in Y'$ with $i < j$, or $y_j \to y_i$ for every $y_i, y_j \in Y'$ with $i < j$.

Let $Y' = \{y_{a(1)}, \ldots, y_{a(t)}\}$, with $a(1) < \ldots < a(t)$, and for each $i \in [t]$ let $y'_i = y_{a(i)}$, so either $y'_1 \to \ldots \to y'_t$, or $y'_t \to \ldots \to y'_1$. Also, for each $i \in [t]$ let $x'_{2i-1} = x_{2a(i)-1}$, and $x'_{2i} = x_{2a(i)}$, and let $X = \{x'_1, \ldots, x'_{2t}\} = \{x_{2i-1}, x_{2i} : y_i \in Y' \}$. Note that $x'_1 \to \ldots \to x'_{2t}$, and $x'_{2i} \to y'_i \to x'_{2i-1}$ for every $i \in [t]$.

We partition $X' \cup Y'$ into blocks of three vertices each, as follows: $D_1 = \{x'_1, x'_2, y'_1\}, D_2 = \{x'_3, x'_4, y'_2\} , \ldots, D_t = \{x'_{2t-1}, x'_{2t}, y'_t\}$, and let $J$ be the complete graph with these $t$ blocks as vertices. Now, define a $16$-colouring $f$ on the edges of $J$ as follows. For each $1 \le i < j \le t$, let $$f(D_iD_j) = 8I[x'_{2i-1} \to y'_j] + 4I[x'_{2i} \to y'_j] + 2I[x'_{2j-1} \to y'_i] + I[x'_{2j} \to y'_i] + 1,$$ where $I[A]$ denotes the indicator function of the event $A$. By Ramsey's Theorem, and our choice of $k$, there exists a complete monochromatic subgraph of $J$ on $s \ge \max\{n, K+1\}$ blocks, in colour $c \in [16]$, say. By renaming the vertices of $G$ if necessary, we may assume that these blocks are $D_1, \ldots, D_s$.

The following claim shows that $c \in \{1,4,13,16\}$.\\

\noindent\ul{Claim}: $(a)$ $x'_1 \to y'_2$ if and only if $x'_2 \to y'_2$, and

\hspace{0.84cm} $(b)$ $x'_3 \to y'_1$ if and only if $x'_4 \to y'_1$.

\begin{proof}
Suppose that $x'_1 \to y'_2$, but $y'_2 \to x'_2$. Since $f$ is monochromatic on $D_1, \ldots, D_s$, it follows that $x'_{2i-1} \to y'_s \to x'_{2i}$ for every $i \in [s-1]$. But now, since $s \ge K + 1$, $T[X' \cup \{y_s\}]$ contains a $K$-structure of Type~1, which contradicts our initial assumption. The proof in the other cases is the same.
\end{proof}

It is now easy to see that $T[X' \cup Y'] = M_I^{(s)}$ for some $I \in \{0,1\}^3$; in fact $I_1 = 1$ if and only if $y'_1 \to y'_2$, $I_2 = 1$ if and only if $x'_1 \to y'_2$, and $I_3 = 1$ if and only if $y'_1 \to x'_3$. Since $s \ge n$, it follows that $M_I^{(n)}$ is a sub-tournament of $T$, and so $M_I^{(n)} \in \P$. Since $n$ was arbitrary, this completes the proof of the lemma.
\end{proof}

It remains to count the number of sub-tournaments of $M_I^{(n)}$ for each $I \in \{0,1\}^3$. For each $n \in \N$, $I \in \{0,1\}^3$ and sufficiently large $m$, let $\L(n,I)$ denote the collection of sub-tournaments of $M_I^{(m)}$ of order $n$. Let $L(n,I) = |\L(n,I)|$ denote the number of distinct tournaments in $\L(n,I)$. The following lemmas cover the various cases. We begin with the simplest case, $I = (1,1,1)$.

\begin{lemma}\label{I=111}
If $I = (1,1,1)$ then $L(n,I) \ge F^*_n$ for every $n \in \N$.
\end{lemma}

\begin{proof}
Let $n \in \N$, $I = (1,1,1)$, $X = \{x_1, \dots, x_{2n}\}$, $Y = \{y_1, \ldots, y_n\}$, and $M = M_I(X,Y)$. For each $i \in [n]$, let $A_i = \{x_{2i-1}, x_{2i}, y_i\}$, and note that $x_{2i-1} \to x_{2i} \to y_i \to x_{2i-1}$. Since $I = (1,1,1)$, we have $A_i \to A_j$ for every $1 \le i < j \le n$, so $M$ consists of an `ordered' set of $n$ cyclic triangles.

We shall map sequences $\{a_1, \ldots, a_m\} \in \{1,3\}^m$ such that $m \in \N$ and $\ds\ds\sum_i a_i = n$, to distinct sub-tournaments of $M$ as follows. If $\omega$ is such a sequence of length $m$, then let $\phi(\omega)$ be induced by the vertices $\{x_{2i} : a_i = 1\} \cup \{x_{2i-1}, x_{2i}, y_i : a_i = 3\}$. To reconstruct $\omega$ from $T = \phi(\omega)$, let $a_1 = 1$ if there is a `top vertex' in $T$ (i.e., a vertex with outdegree $|T| - 1$), and otherwise let $a_1 = 3$ and note that there is a `top triangle' in $T$ (i.e., a triple $\{a,b,c\} \subset V(T)$ with $a \to b \to c \to a$, and $\{a,b,c\} \to d$ for every vertex $d \in V(T) \setminus \{a,b,c\}$). Remove the top vertex/triangle, and repeat, to obtain $a_2$, $a_3$, and so on. There are $F^*_n$ sequences as described, so this proves the result.
\end{proof}

We next consider the case $I_2 \neq I_3$. Given $n \in \N$, $0 \le t \le n/2$, $S \subset [n]$ with $|S| = 2t$, and a $t$-permutation $\sigma$, let $T^*_n(S,\sigma)$ denote the following tournament. Suppose $S = \{a(1), \ldots, a(2t)\}$, with $a(1) < \ldots < a(2t)$. Then $T^*_n(S,\sigma)$ has $n$ vertices, $\{x_1, \ldots, x_n\}$ say, and $x_i \to x_j$ whenever $i < j$, unless $i = a(\ell)$ and $j = a(t + \sigma(\ell))$ for some $\ell \in [t]$, in which case $x_j \to x_i$. Thus $T^*_N(S,\sigma)$ is a transitive tournament, with $t$ independent edges reversed.

\begin{obs}\label{Tstar}
Let $n,t \in \N \cup \{0\}$, and $S = \{a(1), \ldots, a(2t)\} \subset [n]$, where $a(1) < \ldots < a(t) < a(t) + 1 < a(t+1) < \ldots < a(2t)$. Let $\sigma$ be a $t$-permutation, and let $T = T^*_n(S,\sigma)$ be as described above. Then $S$ can be reconstructed from $T$.
\end{obs}

\begin{proof}
Let $n$, $t$, $S$ and $\sigma$ be as described, and let $T$ have vertex set $\{x_1, \ldots, x_n\}$, with $x_i \to x_j$ whenever $i < j$, unless $i = a(\ell)$ and $j = a(t + \sigma(\ell))$ for some $\ell \in [t]$, in which case $x_j \to x_i$. Recall that $d(x)$ denotes the outdegree of $x$, and note that for each $i \in [n]$ we have
\begin{equation*}
d(x_i) =
\begin{cases}
n - i & \text{\emph{if} $i \notin S$,}\\[+1ex]
n - i - 1 & \text{\emph{if} $x_i = a(j)$ and $j \le t$,}\\[+1ex]
n - i + 1 & \text{\emph{if} $x_i = a(j)$ and $j > t$.}
\end{cases}
\end{equation*}
So, $d(x_i) \ge d(x_{i+1})$ for every $i \in [n-1]$, and $x_i \to x_{i+1}$ for every $i \in [n-1]$, since $a(t) + 1 < a(t+1)$. Notice also that at most two vertices can have the same degree, unless $a(t+1) = a(t) + 2$, in which case $x_{a(t)}$, $x_{a(t)+1}$ and $x_{a(t+1)}$ are the only triple all with the same degree.

Now, in order to reconstruct $S$ from $T$, simply order the vertices of $T$ according to their outdegrees, and if two vertices have the same outdegree then order them according to the orientation of the edge between them. More precisely, if $u,v \in V(T)$, then let $u > v$ if $d(u) > d(v)$, or if $d(u) = d(v)$ and $u \to v$. By the comments above, this gives a linear order unless $a(t+1) = a(t) + 2$; in that case, order the three vertices with the same degrees arbitrarily. Let this ordering be $y_1 > \ldots > y_n$, where $\{y_1, \ldots, y_n\} = V(T)$. Now $y_i \to y_j$ if and only if $x_i \to x_j$, so $S = \{y_i : d(y_i) \neq n - i\}$.
\end{proof}

We can now give a lower bound for the case $I_2 \neq I_3$.

\begin{lemma}\label{I1I2}
If $I_2 \neq I_3$ then $L(n,I) \ge 2^{n-2}$ for every $n \in \N$.
\end{lemma}

\begin{proof}
The proofs for the two cases $\{I_2,I_3\} = \{0,1\}$ are almost identical, so for simplicity we assume that $I_2 = 1$ and $I_3 = 0$. Let $n \in \N$, $I_1 \in \{0,1\}$ and $m = n^2 + 2n$. Let $X' = \{x'_1, \dots, x'_{2m}\}$, $Y = \{y_1, \ldots, y_m\}$, and $M' = M_I(X',Y)$. We shall only need to use half of the vertices of $X'$, so for each $i \in [m]$ let $x_i = x'_{2i - 1}$, and let
$$X \; = \; \{x_1, \ldots, x_m\} \; = \; \{x'_i \in X' \; : \; i \equiv 1\textup{ (mod 2)} \}.$$
Note that $y_i \to x_j$ if and only if $i = j$, and let $M = M'[X \cup Y]$. We shall show that $M$ has at least $2^{n-2}$ distinct sub-tournaments.

For each subset $S \subset [n-1]$ of even size, we shall find a distinct sub-tournament $\psi(S)$ of $M$. Let the elements of $S$ be $a(1) < \ldots < a(2t)$, where $0 \le t \le (n - 1)/2$. Let $a(0) = 0$ and $a(2t+1) = n$, and define $b(i) = a(i) - a(i-1) - 1$ for each $i \in [2t+1] \setminus \{t+1\}$, and $b(t+1) = a(t+1) - a(t) \ge 1$. We choose a subset $A \cup B \subset X \cup Y$ on which $M$ is transitive, except for $t$ independent edges which are reversed, as follows. First, let
$$A = \{x_i \in X, \; y_i \in Y \; : \; i = 2jn, \; j \in [t]\}.$$
The vertices of $A$ will be the endpoints of the reversed edges. In order that they correspond to the set $S$, we define a set $B$ as follows. Since our choice of $B$ will depend on the value of $I_1$, we split into two cases.\\[-1.5ex]

\noindent \textbf{Case 1}: $I_1 = 0$.\\[+0.5ex]
\noindent Recall that if $I_1 = 0$, then $y_j \to y_i$ for every $i < j$. Let
\begin{align*}
B \; = \; \{x_i \in X & : \; i = 2jn + \ell, \; j \in [0,t], \; \ell \in [b(j + 1)]\}\\
\cup \;\; \{y_i \in Y & : \; i = (2j-1)n + \ell, \; j \in [t], \; \ell \in [b(2t - j + 2)]\},
\end{align*}
\noindent and let $\,\psi(S) = M[A \cup B]\,$ (see Figure 2). If $S = \emptyset$ then $b(1) = n$, and so $\psi(S)$ is transitive. Otherwise, $b(i) < n$ for every $i \in [2t + 1]$, so $A$ and $B$ are disjoint, and so $|A \cup B| = 2t + \ds\ds\sum_i b(i) = n$.

Now, we claim that $\psi(S) = T^*_n(S',\sigma^{rev})$, where
$$S' = \{a(1), \ldots, a(t), a(t+1)+1, \ldots, a(2t) + 1\} \subset [n],$$ and $\sigma^{rev} = t(t-1) \ldots 21$ is the reverse permutation on $[t]$. Indeed, giving $A \cup B$ the order induced by $x_1 < \ldots < x_m < y_m < \ldots < y_1$, we see that the only `reversed' edges ($u \to v$ but $u > v$) in $\psi(S)$ are those of the form $y_i \to x_i$ with $x_i, y_i \in A$. One can easily check that these correspond to the edges between the $a(\ell)^{th}$ and $(a(t + \sigma^{rev}(\ell)) + 1)^{th}$ elements of $A \cup B$ (in the above ordering) for some $\ell \in [t]$. So $\psi(S) = T^*_n(S',\sigma^{rev})$, as claimed.

So, by Observation~\ref{Tstar} we can reconstruct $S'$ from $\psi(S)$, and it is easy to see that we can reconstruct $S$ from $S'$. Thus $\psi$ is injective. There are $2^{n-2}$ sets $S \subset [n-1]$ of even size, and so this proves the result when $I_1 = 0$.\\[-1.5ex]

\noindent \textbf{Case 2}: $I_1 = 1$.\\[+0.5ex]
\noindent The proof is very similar to that in Case 1. Let
\begin{align*}
B \; = \; \{x_i \in X & : \; i = (2j - 1)n + \ell, \; j \in [t + 1], \; \ell \in [b(j)]\}\\
\cup \;\; \{y_i \in Y & : \; i = 2jn + \ell, \; j \in [t], \; \ell \in [b(t + j + 1)]\},
\end{align*}
and let $\psi(S) = M[A \cup B]$. Now $\psi(S) = T^*_n(S',\sigma^{id})$, where, as before, $S' = \{a(1), \ldots, a(t), a(t+1)+1, \ldots, a(2t) + 1\}$, and $\sigma^{id} = 12 \ldots (t-1)t$ is the identity permutation on $[t]$. Again we can reconstruct $S$ from $\psi(S)$, so this proves the result when $I_1 = 1$.
\end{proof}

\vspace{0.3in}
\[ \unit = 0.47cm \hskip -17\unit
\medline \dl{0}{2}{1}{2} \dl{1}{2}{1}{0} \dl{0}{0}{1}{0} \dl{0}{0}{0}{2}
\thnline \dl{0.1}{1.9}{0.9}{0.1} \dl{0.1}{1}{0.5}{0.1} \dl{0.5}{1.9}{0.9}{1}
\point{12}{-2.6}{$ \medline \dl{0}{2}{1}{2} \dl{1}{2}{1}{0} \dl{0}{0}{1}{0} \dl{0}{0}{0}{2}
\thnline \dl{0.1}{1.9}{0.9}{0.1} \dl{0.1}{1}{0.5}{0.1} \dl{0.5}{1.9}{0.9}{1} $}
\varpt{5000}
\point{0}{-3}{$ \pt{0.5}{-1} \point{-1.6}{-1.2}{$x_{2n}$} \medline \bez{0.5}{-1}{6.5}{1}{12.5}{-1}
\dl{6.3}{0}{6.8}{0.3} \dl{6.3}{0}{6.8}{-0.3} \pt{12.5}{-1} \point{13.2}{-1.4}{$y_{2n}$} $}
\point{0}{-7.3}{$ \medline \dl{0}{2}{1}{2} \dl{1}{2}{1}{0} \dl{0}{0}{1}{0} \dl{0}{0}{0}{2}
\thnline \dl{0.1}{1.9}{0.9}{0.1} \dl{0.1}{1}{0.5}{0.1} \dl{0.5}{1.9}{0.9}{1}
\point{12}{-2.6}{$ \medline \dl{0}{2}{1}{2} \dl{1}{2}{1}{0} \dl{0}{0}{1}{0} \dl{0}{0}{0}{2}
\thnline \dl{0.1}{1.9}{0.9}{0.1} \dl{0.1}{1}{0.5}{0.1} \dl{0.5}{1.9}{0.9}{1} $}
\point{0}{0.5} {$ \medline \dl{5}{0}{8}{0} \dl{8}{0}{7.5}{0.3} \dl{8}{0}{7.5}{-0.3} $}
\varpt{5000}
\point{0}{-3}{$ \pt{0.5}{-1} \point{-1.6}{-1.2}{$x_{4n}$} \medline \bez{0.5}{-1}{6.5}{1}{12.5}{-1}
\dl{6.3}{0}{6.8}{0.3} \dl{6.3}{0}{6.8}{-0.3} \pt{12.5}{-1} \point{13.2}{-1.4}{$y_{4n}$} $} $}
\thkline \dline{6.5}{-12.5}{6.5}{-15.5}{4}
\point{0}{-18}{$ \medline \dl{0}{2}{1}{2} \dl{1}{2}{1}{0} \dl{0}{0}{1}{0} \dl{0}{0}{0}{2}
\thnline \dl{0.1}{1.9}{0.9}{0.1} \dl{0.1}{1}{0.5}{0.1} \dl{0.5}{1.9}{0.9}{1}
\point{12}{-2.6}{$ \medline \dl{0}{2}{1}{2} \dl{1}{2}{1}{0} \dl{0}{0}{1}{0} \dl{0}{0}{0}{2}
\thnline \dl{0.1}{1.9}{0.9}{0.1} \dl{0.1}{1}{0.5}{0.1} \dl{0.5}{1.9}{0.9}{1} $}
\point{0}{0} {$ \medline \dl{5}{0}{8}{0} \dl{8}{0}{7.5}{0.3} \dl{8}{0}{7.5}{-0.3} $}
\varpt{5000}
\point{0}{-3}{$ \pt{0.5}{-1} \point{-1.6}{-1.2}{$x_{2tn}$} \medline \bez{0.5}{-1}{6.5}{1}{12.5}{-1}
\dl{6.3}{0}{6.8}{0.3} \dl{6.3}{0}{6.8}{-0.3} \pt{12.5}{-1} \point{13.2}{-1.4}{$y_{2tn}$} $} $}
\point{0}{-25.3}{$ \medline \dl{0}{2}{1}{2} \dl{1}{2}{1}{0} \dl{0}{0}{1}{0} \dl{0}{0}{0}{2}
\thnline \dl{0.1}{1.9}{0.9}{0.1} \dl{0.1}{1}{0.5}{0.1} \dl{0.5}{1.9}{0.9}{1} $}
\point{18.5}{-20}{$ \medline \dl{0}{2}{1}{2} \dl{1}{2}{1}{0} \dl{0}{0}{1}{0} \dl{0}{0}{0}{2}
\thnline \dl{0.1}{1.9}{0.9}{0.1} \dl{0.1}{1}{0.5}{0.1} \dl{0.5}{1.9}{0.9}{1} \point{1.1}{0.8}{ $= B$} $}
\medline \dl{-3.5}{-9}{-3.5}{-15} \dl{-3.5}{-15}{-3.2}{-14.5} \dl{-3.5}{-15}{-3.8}{-14.5}
\point{19.5}{0} {$ \medline \dl{-3.5}{-9}{-3.5}{-15} \dl{-3.5}{-9}{-3.2}{-9.5} \dl{-3.5}{-9}{-3.8}{-9.5} $}
\point{2}{-27.5}{Figure 2: The set $A \cup B$}
\]
\vspace{0.4in}

We have dealt with the cases $I = (1,1,1)$ and $I_2 \neq I_3$. We next turn to the cases $I = (0,0,0)$ and $I = (0,1,1)$.

\begin{lemma}\label{I1=0}
If $I_1 = 0$ and $I_2 = I_3$, then $L(n,I) \ge 2^{n-3} - 2$ for every $n \in \N$.
\end{lemma}

\begin{proof}
The proofs for the two cases $I = (0,0,0)$ and $I = (0,1,1)$ are almost identical, so we shall only prove the result for $I = (0,0,0)$. Let $n \in \N$ and $m = n^2$. Let $X' = \{x'_1, \dots, x'_{2m}\}$, $Y = \{y_1, \ldots, y_m\}$, and $M' = M_I(X',Y)$. To make the proof easier to follow, we let $x_i = x'_{2i-1}$ for each $i \in [m]$, and let $X = \{x_1, \ldots, x_m\}$. Note that now $x_i \to y_j$ if and only if $i > j$. Let $M = M'[X \cup Y]$.

For each subset $S \subset [n-3] \setminus \{\emptyset, [n-3]\}$, we shall find a distinct sub-tournament $\psi(S)$ of $M$. Let the elements of $S$ be $a(1) < \ldots < a(t)$, where $1 \le t \le n - 4$, and let $a(0) = 0$ and $a(t+1) = n - 2$. For each $i \in [t + 1]$, let $b(i) = a(i) - a(i-1) - 1$, and let
$$A \; = \; \{x_i \in X \; : \; i = (j + 1)n, \; j \in [t]\}.$$
The vertices of $A$ will, as usual, correspond to the elements of $S$. Now let
$$B \; = \; \{y_i \in Y \; : \; i = jn + \ell, \; j \in [t+1], \; \ell \in [b(j)]\},$$
and finally let $$C = \{x_1, y_1, x_m\}.$$
Let $\psi(S) = M[A \cup B \cup C]$. The sets $A$, $B$ and $C$ are disjoint, so $|A \cup B \cup C| = t + \ds\sum b(j) + 3 = n$.

Note that the vertex $y_1$ has outdegree 1 in $\psi(S)$ ($y_1 \to x_1$, but $y_j \to y_1$ and $x_j \to y_1$ if $j > 1$, since $I = (0,0,0)$). We claim that every other vertex in $\psi(S)$ has outdegree at least 2. Indeed, $A \to \{y_1, x_m\}$, so $d(v) \ge 2$ for $v \in A$; $B \to \{x_1,y_1\}$, so $d(v) \ge 2$ for every $v \in B$; $x_1 \to A \cup \{x_m\}$, and $|A| \ge 1$ since $S \neq \emptyset$, so $d(x_1) \ge 2$; and finally, $x_m \to B \cup \{y_1\}$, and $|B| \ge 1$ since $S \neq [n-3]$, so $d(x_m) \ge 2$. Hence $y_1$ is the unique vertex in $\psi(S)$ with degree 1.

So we can identify $y_1$, and since $x_1$ is the unique element of $\{v \in \psi(S) : y_1 \to v\}$, we can also identify $x_1$. But now $A = \{v \in \psi(S) : v \notin \{x_1, y_1\}, x_1 \to v\}$, and moreover $\psi(S)$ is transitive on $A$ and on $B$, and $x_i \to y_j$ if and only if $i > j$. It is now easy to see that we can reconstruct $S$ from $\psi(S)$, and so $\psi$ is injective. This proves the lemma.
\end{proof}

There is only one case left to deal with. Given $n \in \N$, let $C_n$ denote the `cyclic tournament' on $n$ vertices, defined as follows. Let $C_n$ have vertex set $\{x_1, \ldots, x_n\}$, and let $x_i \to x_j$ in $C_n$ if and only if $1 \le j - i < n/2 \pmod n$, or $j - i = n/2$.

\begin{lemma}\label{cycle}
Let $n \in \N$. Then $C_{2n}$ has at least $\ds\frac{2^{n-1}}{n}$ distinct sub-tournaments on $n$ vertices.
\end{lemma}

\begin{proof}
Let $n \in \N$, and let $C = C_{2n}$ be the cyclic tournament defined above, with vertices $\{x_1, \ldots, x_{2n}\}$, and $x_i \to x_j$ if and only if $1 \le j - i \le n - 1 \pmod {2n}$ or $j = i + n$. We shall describe a map $\mu$ from the subsets of $[n-1]$ to sub-tournaments of $C$ on $n$ vertices, such that at most $n$ subsets map to the same tournament, i.e., $|\mu^{-1}(T)| \le n$ for every tournament $T$. Since there are $2^{n-1}$ such subsets, this will suffice to prove the lemma.

For each subset $S \subset [n-1]$, let $A_S = \{x_i : i \in S\}$, $B_S = \{x_{n + i} : i \notin S\}$, and define $\mu(S)$ to be the sub-tournament of $C$ induced by the vertices $A_S \cup B_S \cup \{x_{2n}\}$. Note that $B_S \to x_{2n} \to A_S$, and that if $x_i \in A_S$ and $x_j \in B_S$, then $x_i \to x_j$ if and only if $n + i > j$. Note also that $C[A_S]$ and $C[B_S]$ are transitive tournaments.

Now, we claim that given a tournament $T$ in the image of $\mu$, and the vertex $x \in T$ corresponding to $x_{2n}$ in $C$, we can reconstruct $S$. Indeed, let $A = \{v \in T : x \to v\}$, and $B = \{v \in T : v \to x\}$. Now, for each $u \in A$, let $$s(u) = |\{v \in A : v \to u\}| + |\{w \in B : u \to w\}| + 1.$$ Then $S = \{s(u) : u \in A\}$.

Thus, given any tournament $T$ in the image of $\mu$, there are at most $n$ subsets $S \subset [n-1]$ for which $T = \mu(S)$, and so $|\mu^{-1}(T)| \le n$ for every tournament $T$. By the comments above, this proves the lemma.
\end{proof}

The final case now follows easily.

\begin{lemma}\label{I1=1}
If $I = (1,0,0)$, then $L(n,I) \ge \displaystyle\frac{2^{n-1}}{n}$ for every $n \in \N$.
\end{lemma}

\begin{proof}
Let $n \in \N$, $I = (1,0,0)$, $X = \{x_1, \dots, x_{2n}\}$, $Y = \{y_1, \ldots, y_n\}$, and $M = M_I(X,Y)$. We claim that $C_{2n}$ is a sub-tournament of $M$. By Lemma~\ref{cycle}, this will suffice to prove the result.

Recall that $y_i \to y_j$ if $i < j$, $x_i \to x_j$ if $i < j$, and $x_i \to y_j$ if and only if $i \ge 2j$. Consider the vertices $Z = \{z_1 \ldots, z_{2n}\}$, where $z_i = x_{2i}$ and $z_{n+i} = y_i$ for each $i \in [n]$, and let $\ell \in [2n]$. If $\ell \in [n]$ then $\{w \in Z : z_\ell \to w\} = \{x_{2i} : i > \ell\} \cup \{y_i : i \le \ell\}$, so $z_\ell \to z_j$ if and only if $j \in [\ell+1, n+\ell]$. Similarly, if $\ell \in [n+1,2n]$ then $z_\ell \to z_j$ if and only if $j \in [\ell+1,2n] \cup [1,\ell - n - 1]$. Thus $M[Z] = C_{2n}$, so $C_{2n}$ is a sub-tournament of $M$, as claimed.
\end{proof}

We finish this section by summarising what we have learnt from it. We shall need the following simple observation.

\begin{obs}\label{olarge}
\begin{enumerate}
\item[]
\item[$(i)$] $2^{n-1} - 2\ds{{n-1} \choose 2} - n \, \ge \, F_n^*$ if $n \ge 6$.
\item[$(ii)$] $2^{n-2} > \: 2^{n-3} - 2 \, \ge \, F_n^*$ if $n \ge 6$.
\item[$(iii)$] $\left\lceil \ds\frac{2^{n-1}}{n} \right\rceil \, \ge \, F_n^*$ if $n \neq 4$.
\end{enumerate}
\end{obs}

\begin{proof}
One can check that the inequalities hold for small cases, and that each function $f(n)$ on the left satisfies $f(n+1) \ge f(n) + f(n-2)$.
\end{proof}

Using Lemmas~\ref{T1}, \ref{t2}, \ref{I=111}, \ref{I1I2}, \ref{I1=0} and \ref{I1=1}, and Observation~\ref{olarge}, we reach the following conclusion.

\begin{cor}\label{types}
Let $\P$ be a hereditary property of tournaments, and suppose that $\P$ contains arbitrarily large structures of Type $1$ or $2$. Then
$$|\P_n| \; \ge \; \min \left\{ 2^{n-1} - 2\ds{{n-1} \choose 2} - n, \, F_n^*, \, 2^{n-3} - 2, \, \left\lceil \ds\frac{2^{n-1}}{n} \right\rceil \right\}$$ for every $n \in \N$, and hence $|\P_n| \ge F_n^*$ for every $n \ge 6$.

Moreover, if also $|\P_n| = F_n^*$ for some $n \ge 7$, then $\P$ contains the property $\T = \{T : T \le M_{(1,1,1)}^{(m)}$ for some $m \in \N\}$.
\end{cor}

\begin{proof}
If $\P$ contains arbitrarily large structures of Type~1, then by Lemma~\ref{T1}, $|\P_n| \ge 2^{n-1} - 2\ds{{n-1} \choose 2} - n$. So assume it does not, in which case, $\P$ must contain arbitrarily large structures of Type~2, and so, by Lemma~\ref{t2}, $M_I^{(n)} \in \P$ for some $I \in \{0,1\}^3$, and every $n \in \N$. Now, by Lemmas~\ref{I=111}, \ref{I1I2}, \ref{I1=0} and \ref{I1=1}, we have
$$|\P_n| \ge \min \left\{ F_n^*, 2^{n-2}, 2^{n-3} - 2, \left\lceil \ds\frac{2^{n-1}}{n} \right\rceil \right\}$$
for every $n \in \N$, and the first statement follows. By Observation~\ref{olarge}, it follows that $|\P_n| \ge F_n^*$ for every $6 \le n \in \N$.

Now suppose that $|\P_n| = F_n^*$ for some $n \ge 7$. Since $2^{n-1} - 2\ds{{n-1} \choose 2} - n > F_n^*$ for $n \ge 7$, it follows that $\P$ does not contain arbitrarily large structures of Type~1, and since $\min\left\{ 2^{n-3} - 2, \, 2^{n-2}, \, \left\lceil \ds\frac{2^{n-1}}{n} \right\rceil \right\} > F_n^*$ if $n \ge 7$, it follows that $\P$ does not contain the set $\{M_I^{(n)} : n \in \N\}$ for any $I \in \{0,1\}^3 \setminus \{(1,1,1)\}$. Thus $M_{(1,1,1)}^{(n)} \in \P$ for every $n \in \N$, and so $\{T : T \le M_{(1,1,1)}^{(m)}$ for some $m \in \N\} \subset \P$, as claimed.
\end{proof}

\section{Polynomial speed}\label{tpolysec}

The results of the previous two sections imply that if $\P$ is a hereditary property of tournaments, and $B(\P) = \infty$, then $|\P_n|$ grows at least exponentially as $n \to \infty$. By Lemma~\ref{Bbdd}, we know that if $B(\P) < \infty$, then $|\P_n|$ is bounded above by a polynomial. In this section we shall prove Theorem~\ref{tpolythm}, which considerably extends Lemma~\ref{Bbdd}. In other words, we shall show that if $B(\P) < \infty$, then $|\P_n| = \Theta(n^k)$, where
\begin{eqnarray*}
k(\P) & = & \sup\{ \ell : \forall \, m \in \N, \: \exists \, T \in \P\textup{ such that the }(\ell + 1)^{st}\textup{ largest}\\
&& \textup{\hspace{1.3cm} homogeneous block in $T$ has at least $m$ elements}\}.
\end{eqnarray*}
We shall also give some idea of how this result might be further improved. The method is very similar to that of Section 4 of~\cite{order}, and the reader may wish to compare the results obtained here with those obtained in that paper for ordered graphs.

We begin by recalling that a pair $\{u,v\}$ of vertices of a tournament $T$ are homogeneous (and we write $u \sim v$) if $u = v$, or $u \curvearrowright v$, or $v \curvearrowright u$, and that the homogeneous blocks of $T$ are the equivalence classes of the relation $\sim$. The \emph{homogeneous block sequence} of $T$ is the sequence $(t_1,t_2,\ldots,t_m)$, where $t_1,  t_2, \ldots, t_m \in \N$ are the orders of the homogeneous blocks of $T$, and $t_1 \ge t_2 \ge \ldots \ge t_m$. Note that this sequence is uniquely determined by $T$. Recall also that $B(T)$ denotes the number of homogeneous blocks of $T$, so $B(T) = m$. We may also embed the homogeneous block sequence of $T$ into the space of infinite sequences of non-negative integers in a natural way, in which case $B(T) = \min\{k \in \N : t_{k+1} = 0\}$.

Now, let $T$ be a tournament, and $B_1, \ldots, B_m$ be the homogeneous blocks of $T$. Define $T(B_1, \ldots, B_m)$ to be the labelled (or ordered) tournament $U$, with vertex set $[m]$, and in which $i \to j$ in $U$ if and only if $b_i \to b_j$ in $T$ for some (and so every) $b_i \in B_i$ and $b_j \in B_j$.

Let $\P$ be a hereditary property of tournaments, and suppose that $B(\P) < \infty$, so there exists a $K \in \N$ such that $B(T) \le K + 1$ for every $T \in \P$. Thus $t_m = 0$ if $m \ge K + 2$, so $\sum_{i=k+2}^\infty t_i$ is bounded for some $k \in \N$ with $k \le K$. The following lemma shows that in this case, $|\P_n| = O(n^k)$.

\begin{lemma}\label{tblocks}
Let $\P$ be a hereditary property of tournaments, and let $k, M \ge 0$ be integers. Suppose that for every $T \in \P$, the homogeneous block sequence of $T$ satisfies $\ds\sum_{i = k+2}^\infty t_i \le M$. Then $|\P_n| = O(n^k)$.
\end{lemma}

\begin{proof}
Let $\P$ be a hereditary property of tournaments, let $k,M \ge 0$ be integers, and suppose that $t_{k+2} + t_{k+3} + \ldots \le M$ for every $T \in \P$. We shall give an upper bound on the number of tournaments of order $n$ in the property.

Indeed, every tournament $T \in \P_n$ is determined by a sequence $S = (a_1, \ldots, a_m)$ of positive integers, satisfying $1 \le m \le k + M + 1$, $\sum_{i=1}^m a_i = n$, and $\sum_{i \in I} a_i \ge n - M$ for some set $I \subset [m]$ with $|I| \le k + 1$; and an ordered tournament $U$ on $m$ vertices. To see this, let $T \in \P_n$ have homogeneous blocks $B_1, \ldots, B_m$, let $a_i = |B_i|$ for each $i \in [m]$, and let $U = T(B_1, \ldots, B_m)$. Now, $1 \le m \le k + M + 1$, since $\sum_{i = k+2}^\infty t_i \le M$; $\sum_{i=1}^m a_i = n$ since $|T| = n$; and $\sum_{i \in I} a_i \ge n - M$ if $I = \{i : B_i$ is one of the largest $k+1$ homogeneous blocks of $T\}$. Thus $S = (a_1, \ldots, a_m)$ and $U$ satisfy the conditions above. It is clear that $T$ can be reconstructed from $S$ and $U$.

It remains to count the number of such pairs $(S,U)$. The number of sequences $S$ is at most
\begin{align*}
& \ds\ds\sum_{m=1}^{k+1} {{n-1} \choose {m-1}} + \ds\sum_{m=k+2}^{k+M+1} {m \choose {k+1}} {{n-M+k} \choose k} {{M + m - 1} \choose {m-1}} \; = \; O(n^k),
\end{align*}
and the number of ordered tournaments on $m$ vertices is just a constant, so this proves the result.
\end{proof}

Theorem~\ref{tpolythm} says that in fact, if $k$ is taken to be minimal in Lemma~\ref{tblocks}, then $|\P_n| = \Theta(n^k)$. The next lemma provides the required lower bound.

\begin{lemma}\label{geton-2}
Let $\P$ be a hereditary property of tournaments, let $k \in \N$, and suppose that there are tournaments $T \in \P$ such that $t_{k+1}$, the size of the $(k+1)^{st}$ largest homogeneous block in $T$, is arbitrarily large. Then
$$|\P_n| \ge \frac{1}{(k+1)!}{{n - 2(k+1)^3} \choose k} = \frac{n^k}{k!(k+1)!} + O(n^{k-1})$$ as $n \to \infty$. Moreover, if $k = 1$, then $|\P_n| \ge n - 2$ for every $n \in \N$.
\end{lemma}

\begin{proof}
Let $\P$ be a hereditary property of tournaments, let $n,k \in \N$, and let $T \in \P$ have $k+1$ homogeneous blocks of order at least $n$. We shall construct a sub-tournament $U$ of $T$ with at least $$\ds\frac{1}{(k+1)!}{{n - 2(k+1)^3} \choose k}$$ distinct ordered sub-tournaments. The idea is simply that $U$ should also have $k+1$ large homogeneous blocks, and at most $3{{k+1} \choose 2}$ other vertices.

Let $B_1, \ldots, B_{k+1}$ be homogeneous blocks of $T$, each of order at least $n$. Let $A_0 = B_1 \cup \ldots \cup B_{k+1}$. We shall inductively define a sequence of sets $A_0 \subset A_1 \subset \ldots \subset A_t$, for some $t \in \left[ 0,{{k+1} \choose 2} \right]$, such that $|A_{i+1}| \le |A_i| + 3$ for each $i \in [1,t-1]$, and so that the sets $\{B_i : i \in [k+1]\}$ are all in different homogeneous blocks of $U = T[A_t]$.

Let $i \in [0,k-1]$, suppose we have already defined the sets $A_0 \subset \ldots \subset A_i$, and let $U_i = T[A_i]$. If the sets $\{B_i : i \in [k+1]\}$ are all in different homogeneous blocks of $U_i$, then we are done with $t = i$ and $U = U_i$. So suppose that there exist $j,\ell \in [k]$ (with $j \neq \ell$) such that $B_j$ and $B_\ell$ are in the same homogeneous block of $U_i$. We shall find a set $A_{i+1}$ as required, such that $B_j$ and $B_\ell$ are in different homogeneous blocks of $T[A_{i+1}]$. Since $B_j$ and $B_\ell$ are distinct homogeneous blocks of $T$, either $B_j \to B_\ell$ or $B_\ell \to B_j$. Without loss of generality, assume that $B_j \to B_\ell$.\\[-1.4ex]

\noindent \textbf{Case 1}: there exists a vertex $u \in T$ such that $B_\ell \to u \to B_j$.\\[+0.6ex]
Let $A_{i+1} = A_i \cup \{u\}$. We claim that $B_j$ and $B_\ell$ are in different homogeneous blocks of $T[A_{i+1}]$. Indeed, let $b \in B_j$ and $b' \in B_\ell$, and suppose that $b \sim b'$, with homogeneous path $C$. We know that $b \to b' \to u \to b$, so $u \notin C$, but now $b' \to u \to b$ is a contradiction. Hence $B_j$ and $B_\ell$ are indeed in different homogeneous blocks of $T[A_{i+1}]$, as required.

So suppose there is no such vertex $u$ with $B_\ell \to u \to B_j$. Let $K = \{v \in T : v \to B_j \cup B_\ell\}$, $L = \{v \in T : B_j \cup B_\ell \to v\}$, and $M = \{v \in T : B_j \to v \to B_\ell\}$. We have $K \cup L \cup M = V(T) \setminus (B_j \cup B_\ell)$.\\[-1.4ex]

\noindent \textbf{Case 2}: $u \to v$ for some $u \in M$ and $v \in K$.\\[0.6ex]
Let $A_{i+1} = A_i \cup \{u,v\}$, and suppose that $b \sim b'$ in $T[A_{i+1}]$, with $b$ and $b'$ as above. Since $b \to u \to b'$, we must have $u \in C$. But now $u \to v \to b$ is a contradiction, so $B_j$ and $B_\ell$ are indeed in different homogeneous blocks of $T[A_{i+1}]$.\\[-1.4ex]

\noindent \textbf{Case 3}: $u \to v$ for some $u \in L$ and $v \in M$.\\[0.6ex]
Let $A_{i+1} = A_i \cup \{u,v\}$. As in Case 2, $B_j$ and $B_\ell$ are in different homogeneous blocks of $T[A_{i+1}]$.\\[-1.4ex]

\noindent \textbf{Case 4}: $u \to v \to w \to u$ for some $u, v, w \in M$.\\[0.6ex]
Let $A_{i+1} = A_i \cup \{u,v,w\}$, and suppose that $b \sim b'$ in $T[A_{i+1}]$, with $b$ and $b'$ as above. Since $b \to u \to b'$, we must have $u \in C$, and similarly, $v,w \in C$. But now $u \to v \to w \to u$ is a contradiction, so $B_j$ and $B_\ell$ are once again in different homogeneous blocks of $T[A_{i+1}]$.\\[-1.4ex]

So suppose that none of the above four cases hold. Then $K \to M \to L$, and $T[M]$ is transitive. It is now easy to see that all the vertices of $B_j \cup B_\ell \cup M$ are in the same homogeneous block of $T$, and this is a contradiction.

We have shown that we can construct sets $A_0 \subset \ldots \subset A_t$ with $|A_{i+1}| \le |A_i| + 3$ for each $i \in [1,t-1]$. Now, the sequence $(A_0, \ldots, A_t) $ cannot continue any further than $t = {{k+1} \choose 2}$, since if $B_j$ and $B_\ell$ are in different homogeneous blocks of $U_i$ (for some $i \in [0,t-1]$ and $j,\ell \in [k+1]$), then they are in different homogeneous blocks of $U_{i+1}$. Since each step of the process described above separates $B_j$ and $B_\ell$ for at least one pair $j,\ell \in [k+1]$, after ${{k+1} \choose 2}$ steps all $k+1$ sets $B_i$ must be in different homogeneous blocks of $U = T[A_t]$.

Now, $U$ has $k+1$ homogeneous blocks of size at least $n$, and at most $K = 3{{k+1} \choose 2}$ other vertices, since $|A_{i+1}| \le |A_i| + 3$ for each $i \in [0,t-1]$. Consider the sub-tournaments of $U$ of order $n$ which include all the vertices of $A_t \setminus A_0$, and $a_i$ vertices from $B_i$ (for each $i \in [k+1]$), where $a_i \ge K + 1$ for each $i \in [k+1]$, and $a_1 \ge \ldots \ge a_{k+1}$. These sub-tournaments are all distinct, since they have different homogeneous block sequences. There are exactly $\ds{{n - (k+1)(K+1) - K + k} \choose k}$ sequences of integers $(a_1, \ldots, a_{k+1})$, with $a_i \ge K + 1$ for each $i \in [k+1]$, and $\ds\sum a_i = n - K$, and at least $\left( \ds\frac{1}{(k+1)!} \right)^{th}$ of these have $a_1 \ge \ldots \ge a_{k+1}$. Therefore there are at least this many distinct sub-tournaments of $U$ of order $n$, and each of these is in $\P$. Finally, note that $(k+1)(K+1) - K + k < 2(k+1)^3$, so $$|\P_n| \ge \frac{1}{(k+1)!}{{n - 2(k+1)^3} \choose k} = \frac{n^k}{k!(k+1)!} + O(n^{k-1}),$$ as required.

To prove the second part of the lemma, let $k = 1$, $n \in \N$, and repeat the argument above to obtain the tournament $U \in \P$, with two homogeneous blocks of order $n$ and at most three other vertices. Now, the four cases in the proof above correspond exactly to the four tournaments $G_1^{(n)}$, $G_2^{(n)}$, $G_3^{(n)}$ and $G_4^{(n)}$ defined in \cite{posets}, so $\P$ must contain one of these tournaments. By a simple counting argument, we have $|\P_n| \ge n - 2$ for every $n \in \N$.
\end{proof}

\begin{rmk}
The constant $\ds\frac{1}{k!(k+1)!}$ in Lemma~\ref{geton-2} is not best possible. In fact, with a little more care one can replace it with $$\left( k!\ds\max_T(|\operatorname{Aut}(T)|) \right)^{-1},$$ where $\operatorname{Aut}(T)$ denotes the automorphism group of $T$, and the maximum is taken over all tournaments on $k+1$ vertices. Consider the following sequence of tournaments: $T_1$ is a cyclic triangle, and for each $\ell \in \N$, $T_{\ell + 1}$ is formed by taking three copies $U$, $V$ and $W$ of $T_\ell$, and letting $U \to V \to W \to U$ in $T_{\ell + 1}$. The automorphism group of $T_\ell$ has size $3^{(k-1)/2}$, where $k = 3^\ell$ is the number of vertices of $T_\ell$, and this was shown to be the largest possible order of the automorphism group of a tournament in 1970 by Moon~\cite{Moon}. Moon's result (together with the argument above) implies that the bound in Lemma~\ref{geton-2} could be improved to $$|\P_n| \; \ge \; \frac{n^k}{k! \: 3^{(k-1)/2}} + O(n^{k-1}),$$ and this constant is in fact best possible.
\end{rmk}

Theorem~\ref{tpolythm} follows easily from Lemmas~\ref{tblocks} and \ref{geton-2}.

\begin{proof}[Proof of Theorem~\ref{tpolythm}]
Let $k = k(\P)$, and suppose $k < \infty$. By the definition of $k(\P)$, $\sum_{i = k+2}^\infty t_i \le M$ for some $M \in \N$, and there are tournaments $T \in \P$ such that $t_{k+1}$, the size of the $(k+1)^{st}$ largest homogeneous block in $T$, is arbitrarily large. Thus $|\P_n| = O(n^k)$ by Lemma~\ref{tblocks}, and $|\P_n| = \Omega(n^k)$ by Lemma~\ref{geton-2}.
\end{proof}

\section{Proof of Theorem~\ref{tourns}}\label{tproof}

Theorem~\ref{tourns} now follows easily from Lemma~\ref{tournkey}, the results of Section~\ref{T1T2}, and Theorem~\ref{tpolythm}. The only remaining ingredient is the following lemma, which covers the case $n \le 5$. We shall only give a sketch of the (easy, but tedious) details of the proof.

\begin{lemma}\label{osmall}
Let $\P$ be a hereditary property of tournaments.
\begin{enumerate}
\item[$(a)$] If $|\P_n| \ge 2$ for some $n \in \N$, then $|\P_n| \ge F_n^*$ for $n = 1,2,3$.
\item[$(b)$] If $\P$ contains a $3$-structure of Type~1, then $|\P_5| \ge F_5^* = 4$.
\item[$(c)$] If $M_I^{(3)} \in \P$, with $I_1 = 0$, then $|\P_5| \ge F_5^* = 4$.
\item[$(d)$] If $\P = \{T : T \le C_n$ for some $n \in \N\}$, then $|\P_4| = 2 < 3 = F_4^*$.
\end{enumerate}
\end{lemma}

\begin{proof}
For part $(a)$, note that if $|\P_n| \ge 2$ for some $n \in \N$, then $\P$ contains both tournaments on three vertices (the transitive tournament and the cyclic triangle). Since $|F_1^*| = |F_2^*| = 1$ and $|F_3^*| = 2$, the result is immediate. For part $(b)$, observe that a 3-structure of Type~1 contains all four tournaments in $\T_5$, so if $\P$ contains such a structure, then $\T_5 \subset \P_5$.  Similarly, for part $(c)$ observe that if $I_1 = 0$, $M_I^{(3)}$ contains a tournament on five vertices with a 5-cycle, both tournaments with a 4-cycle (one with a vertex `above' the 4-cycle, and one with a vertex `below' it), and the transitive tournament. Finally, for part $(d)$ note that $C_n$ does not contain a cyclic triangle for $n \ge 4$, so $\P_4$ contains only the transitive tournament and the 4-cycle.
\end{proof}

Finally, we are ready to prove Theorem~\ref{tourns}.

\begin{proof}[Proof of Theorem~\ref{tourns}]
Let $\P$ be a hereditary property of tournaments, and suppose first that $B(\P) < \infty$. Then $k = k(\P) \le B(\P)$, so $|\P_n| = \Theta(n^k)$ as $n \to \infty$ by Theorem~\ref{tpolythm}. So assume that $B(\P) = \infty$. By Lemma~\ref{tournkey}, $\P$ contains arbitrarily large structures of Type~1 or~2 and so by Corollary~\ref{types},
$$|\P_n| \; \ge \; \min \left\{ 2^{n-1} - 2\ds{{n-1} \choose 2} - n, \, F_n^*, \, 2^{n-3} - 2, \, \left\lceil \ds\frac{2^{n-1}}{n} \right\rceil \right\}$$ for every $n \in \N$, and hence $|\P_n| \ge F_n^*$ if $n \ge 6$.

It only remains to show that $|\P_n| \ge F_n^*$ for $n \le 5$ and $n \neq 4$. For $n \in \{1,2,3\}$ this follows trivially by Lemma~\ref{osmall} $(a)$, but for $n = 5$ we must do a tiny bit of work. Recall that since $B(\P) = \infty$, either $\P$ contains arbitrarily large structures of Type 1, or $\P$ contains the tournament $M_I^{(n)}$ for some $I \in \{0,1\}^3$ and every $n \in \N$ (see Lemmas~\ref{tournkey} and \ref{t2}). In the former case, we have $|\P_5| \ge F_5^*$ by Lemma~\ref{osmall} $(b)$. In the latter case, we have $|\P_5| \ge F_5^*$ if $I = (1,1,1)$ (by Lemma~\ref{I=111}), if $I_2 \neq I_3$ (by Lemma~\ref{I1I2}), and if $I = (1,0,0)$ (by Lemma~\ref{I1=1} and Observation~\ref{olarge}). But $I_1 = 0$ in the remaining cases, and so the result follows by Lemma~\ref{osmall} $(c)$. Finally, Lemma~\ref{osmall} $(d)$ shows that the hereditary property $\P = \{T : T \le C_n$ for some $n \in \N\}$, which by Lemma~\ref{cycle} has speed at least $2^{n-1}/n$, satisfies $|\P_4| < F_4^*$.
\end{proof}

\section{Further problems}\label{tprobs}

Research into hereditary properties of tournaments is still at an early stage, and we have many more questions than results. We present here a selection of problems and conjectures; we begin with a Stanley-Wilf Conjecture for tournaments.

\begin{conj}
There is a jump from exponential to factorial speed for hereditary properties of tournaments. More precisely, there exists a constant $\alpha > 0$ and a function $F(n) = n^{\alpha n + o(n)}$, such that, for any hereditary property of tournaments $\P$, either\\[-2.5ex]
\begin{enumerate}
\item[$(a)$] $|\P_n| \le c^n$ for every $n \in \N$, for some constant $c = c(\P)$, or\\[-1.5ex]
\item[$(b)$] $|\P_n| \ge F(n)$ for every $n \in \N$.
\end{enumerate}
\end{conj}

We would also like to know which exponential speeds are possible.

\begin{qu}\label{tctothen}
Let $\P$ be a hereditary property of tournaments, and suppose that $|\P_n| < c^n$ for some $c \in \RR$ and every $n \in \N$. Does $\ds\lim_{n \to \infty} \left( |\P_n|^{1/n} \right)$ necessarily exist?
\end{qu}

\begin{prob}\label{tbases}
Let $\X = \{c \in \RR :$ there is a hereditary property of tournaments $\P$ with $\ds\lim_{n \to \infty} \left( |\P_n|^{1/n} \right) = c\}$. Determine the set $\X$.
\end{prob}

Theorem~\ref{tourns} implies that $\X \cap (0, c_3) = \emptyset$, where $c_3 \simeq 1.47$ is the largest real root of the polynomial $x^3 = x^2 + 1$. But what happens above $c_3$? Consider the following generalization of the tournament $M_{(1,1,1)}^{(n)}$. Let $k,n \in \N$ with $k \ge 3$, and let $M(k,n)$ be the tournament with vertex set $\{x_1, \ldots, x_{kn}\}$, in which $x_i \to x_j$ if $i < j$, unless $i + k - 1 = j \equiv 0 \pmod k$. Note that $M(3,n) = M_{(1,1,1)}^{(n)}$. Now let $\P^{(k)} = \{T : T \le M(k,n)$ for some $n \in \N\}$, and observe that $|\P^{(k)}_n| = c_k^{(1+o(1))n}$, where $c_k$ is the largest real root of the polynomial $x^k = x^{k-1} + x^{k-3} + x^{k-4} + \ldots + 1$. Note also that $c_k \to c'$ as $k \to \infty$, where $c'$ is the largest real root of the polynomial $x^4 = x^3 + x^2 + 1$.

We conjecture, along the lines of Theorem 1 of~\cite{order}, that these are the only bases in the range $[0,c']$.

\begin{conj}
Let $\X$ be as defined in Problem~$\ref{tbases}$, let $c' \simeq 1.755$ be the largest real root of the polynomial $x^4 = x^3 + x^2 + 1$, and let $A = \{x : x$ is the largest real root of the polynomial $x^{k+1} = x^k + x^{k-2} + x^{k-3} + \ldots + 1$ for some $2 \le k \in \N\}.$ Then $\X \cap [0,c'] = \{0,c'\} \cup A.$
\end{conj}

As in \cite{order}, we can prove a special case of Question~\ref{tctothen}. A tournament $T$ is said to be \emph{strongly connected} if every ordered pair of vertices $u$ and $v$ are connected by a path from $u$ to $v$. Equivalently, the vertex set of $T$ cannot be partitioned into two non-empty sets $A$ and $B$ such that $A \to B$.

\begin{thm}\label{tfekete}
Let $T_1, T_2, \ldots$ be a sequence of tournaments, and suppose that every $T_i$ is strongly connected. Let $\P = \{T : T$ is a tournament, and $T_i \not\le T$ for every $i \in \N\}$. Then either $\ds \lim_{n \to \infty} \left( |\P_n|^{1/n} \right)$ exists, or $\ds \liminf_{n \to \infty} \left( |\P_n|^{1/n} \right) = \infty$.
\end{thm}

\begin{proof}
The proof is essentially the same as that of Theorem 27 in~\cite{order}. We claim that for every pair of integers $m,n$, $$|\P_{m+n}| \ge |\P_m| \cdot |\P_n|.$$ To see this, let $G_1 \in \P_m$ and $G_2 \in \P_n$, and let $(G_1,G_2)$ denote the tournament on $m + n$ vertices formed by taking disjoint copies of $G_1$ and $G_2$, and orienting all cross-edges from $V(G_1)$ to $V(G_2)$. Then $T_i \notin (G_1,G_2)$ for every $i$, so $(G_1,G_2) \in \P_{m+n}$, and moreover $G_1$ and $G_2$ can be reconstructed from $(G_1,G_2)$, so the claim follows.

Now, Fekete's Lemma~\cite{Fek} states that if $a_1, a_2, \ldots \in \RR$ satisfy $a_m + a_n \ge a_{m+n}$ for all $m,n \ge 1$, then $\ds\lim_{n \rightarrow \infty} \displaystyle\frac{a_n}{n}$ exists and is in $[-\infty,\infty)$. Applying this lemma to the sequence $-\log(|\P_n|)$ gives the result.
\end{proof}

The proof of Theorem 28 of \cite{order} can also be adapted to hereditary properties of tournaments, to produce many properties with different exponential speeds, but we spare the reader the details. Of perhaps more interest is whether our results from this paper can be used to prove a jump from polynomial to exponential speed for hereditary properties of (unlabelled) oriented and directed graphs. We therefore finish with the following conjecture.

\begin{conj}
Let $\P$ be a hereditary property of oriented graphs. Then either
\begin{enumerate}
\item[$(a)$] $|\P_n| = \Theta(n^k)$ for some $k \in \N$, or\\[-1.5ex]
\item[$(b)$] $|\P_n| \ge F^*_n$ for every $4 \neq n \in \N$.\\[-2ex]
\end{enumerate}
\end{conj}

\end{document}